\documentclass[12pt]{article}
\usepackage{amsfonts,amsmath,amssymb}
\usepackage{enumerate}
\usepackage{theorem}
\numberwithin{equation}{section}

\newtheorem{theorem}{Theorem}[section]
\newtheorem{lemma}[theorem]{Lemma}
\newtheorem{proposition}[theorem]{Proposition}

\newtheorem{example}[theorem]{Example}

\newtheorem{definition}{Definition}[section]

\newtheorem{remark}[theorem]{Remark}

\newcommand{\cl}[1]{\mathcal{#1}} 
\newcommand{\bb}[1]{\mathbb{#1}}

\newcommand{\nor}[1]{\left\Vert #1\right\Vert}

\begin{document}

\title{MORITA TYPE EQUIVALENCES AND REFLEXIVE ALGEBRAS}
\author{
G.K. ELEFTHERAKIS\\
 Department of Mathematics,\\ University of
Athens, \\ Panepistimioupolis 157 84\\ Athens Greece\\
e-mail: gelefth@math.uoa.gr
}

\date{}

\maketitle

\begin{abstract}
Two unital dual operator algebras $\cl{A}, \cl{B}$ are called $\Delta $-equivalent 
if there exists an equivalence  functor
$\cl{F}:\; _{\cl{A}}\mathfrak{M}\rightarrow \; _{\cl{B}}\mathfrak{M}$
which ``extends" to a $*-$functor implementing an equivalence between the categories
$_{\cl{A}}\mathfrak{DM}$ and $_{\cl{B}}\mathfrak{DM}.$ Here $_{\cl{A}}\mathfrak{M}$
denotes the category of normal representations of $\cl{A}$ and $_{\cl{A}}\mathfrak{DM}$
 denotes the category with the same objects as $_{\cl{A}}\mathfrak{M}$ and $\Delta (\cl{A})$-module
maps as morphisms ($\Delta (\cl{A})=\cl{A}\cap \cl{A}^*$).
We prove that any such functor maps completely isometric
representations to completely isometric representations, ``respects" the
lattices of the algebras and maps reflexive algebras to reflexive algebras. We 
present applications to the class of CSL algebras.
\end{abstract}



\section{Introduction.}

The Morita equivalence of rings has been extended to many settings. In particular, Rieffel 
developed the appropriate notion and theory for $C^*$-algebras and $W^*$-algebras 
\cite{rif}, \cite{rif2}. Blecher, Muhly and Paulsen generalized Rieffel's $C^*$-algebraic 
version to nonselfadjoint operator algebras \cite{bmp}. In \cite{ele}, \cite{ele2} we obtained
 a generalization of Rieffel's concept 
of Morita equivalence of $W^*$-algebras to the class of (not necessarily selfadjoint) 
unital dual operator algebras. The purpose of the present work is to apply this theory to the case of reflexive 
algebras, especially to CSL algebras. As we know this is the first connection 
between Morita theory and nonselfadjoint reflexive algebras in literature.

We say that two unital operator algebras are $\Delta $-equivalent if there is an equivalence functor 
between their categories of normal representations which not only preserves intertwiners of 
representations of the algebras, but also preserves intertwiners of their restrictions to the diagonals 
(see Definition \ref{1.4.d}). This reduces to Rieffel's definition in the selfadjoint case. 
In section 2 we study properties of $\Delta$-equivalence functors. 
We show that every such functor maps completely isometric representations to completely isometric 
ones. Furthermore, every equivalence functor sends the invariant projection lattice
of a representation  onto the invariant projection lattice
of the corresponding representation,
and maps reflexive algebras to reflexive algebras.

In section 3 we give examples of $\Delta$-equivalent algebras and of $\Delta$-inequivalent algebras.
We show that two CSL algebras $\cl{A, B}$ are $\Delta$-equivalent if  and only if 
there exists a ternary ring of operators $\cl{M}$ such that $\cl{A}=[\cl{M}^*\cl{B}\cl{M}]^{-w^*}\;\; \text{and}\;\;
\cl{B}=[\cl{M}\cl{A}\cl{M}^*]^{-w^*}$ (see Definition \ref{1.1.d}.)
 It follows from \cite{ele} that two separably acting CSL algebras with either continuous or totally atomic invariant projection
lattices are $\Delta$-equivalent if  and only if  they have isomorphic lattices.
We show however that isomorphism of the invariant projection lattices 
does not always imply $\Delta$-equivalence of the algebras,
even in the case of multiplicity free nests and isomorphic diagonals.
Nevertheless, if  there exists an order isomorphism between two nests, preserving the dimensions
of intervals (equivalently, if the nests  are {\em similar}), there always exists an
equivalence functor between the categories of normal, completely contractive representations
of the nest algebras,
which is normal and completely isometric.

We inform the reader that in another paper \cite{elepaul}, written after the 
research reported in this paper was completed, we proved jointly 
with Vern Paulsen that two unital dual operator algebras $\cl{A, B}$ are 
$\Delta $-equivalent if and only if they are stably isomorphic, 
i.e. if there exists a Hilbert space $H$ such that the algebras 
$\cl{A} \overline{\otimes }B(H)$ and $\cl{B}\overline{\otimes }B(H)$ 
(where $\overline{\otimes }$ denotes the normal spatial tensor product \cite{bm}) 
 are isomorphic as dual operator algebras.

We present some symbols used below. If $\cl{A}$ is an operator algebra \cite{bm}, \cite{paul},
 we denote its diagonal $\cl{A}\cap \cl{A}^*$ by $\Delta (\cl{A}).$ 
The symbol $[\cl{S}]$ denotes the linear span of
$\cl{S}.$ The commutant of a set $\cl{L}$ of bounded operators on a 
Hilbert space $H$ is denoted $\cl{L}^\prime.$ If $\cl{U}$ 
is a linear space and $n,m\in \mathbb{N}$ we denote by $M_{n,m}(\cl{U})$ the space of
 $n\times m$ matrices with entries from $\cl{U}$ and by $M_{n}(\cl{U})$ the space 
$M_{n,n}(\cl{U}).$ If $\cl{U}, \cl{V}$ are linear spaces, $\alpha $ is  a linear map 
from $\cl{U}$ to $\cl{V}$ and $n,m\in \mathbb{N}$ we denote the linear map\vspace*{-1ex} 
$$M_{n,m}(\cl{U})\rightarrow  M_{n,m}(\cl{V}): (A_{ij})_{i,j}\rightarrow 
(\alpha (A_{ij}))_{i,j}$$ 
again by $\alpha .$ If $\cl{U}$ is a subspace of 
$B(H,K)$ for $H,K$ Hilbert spaces we equip $M_{n,m}(\cl{U}) $, $n,m\in \mathbb{N}$ 
with the norm inherited from the embedding $M_{n,m}(\cl{U})\subset B(H^m,K^n).$ If 
$(\cl{X}, \|\cdot\|)$ is a normed space we denote by $Ball(\cl{X})$ the unit ball of
 $\cl{X}:$ $\{X\in \cl{X}: \|X\| \leq 1 \}.$ If $x_1,...,x_n$ are in a vector space 
$\cl{V},$ we write $(x_1,...,x_n)^t$ for the column vector in $M_{n,1}(\cl{V}). $

A $C^*$ algebra which is a dual Banach space is called a \textbf{$W^*$ algebra}. 
A \textbf{dual operator algebra} is an abstract operator algebra which is the operator dual of an 
operator space. Every $W^*$ algebra is a dual operator algebra. For 
every dual operator algebra $\cl{A}$ there exists a Hilbert space $H_0$ and an 
algebraic homomorphism $\alpha _0: \cl{A}\rightarrow B(H_0)$ which is a complete 
isometry and a $w^*$-continuous map \cite{bm}.

 A set of projections on a Hilbert space will be called a \textbf{lattice} if contains the zero and identity operators and is closed under arbitrary suprema and infima. If $\cl{A}$ is a subalgebra of $B(H)$ for some Hilbert space $H,$ the set
$$\mathrm{Lat}(\cl{A})=\{L\in pr(B(H)): L^\bot \cl{A} L=0\}$$ is a lattice. Dually if $\cl{L}$ is a lattice the space $$\mathrm{Alg}(\cl{L})=\{A\in B(H): L^\bot AL=0 \; \forall \; L\in \cl{L}\}$$ is an algebra. The \textbf{reflexive hull} of a unital algebra $\cl{A}$ is the algebra $$\mathrm{Ref}(\cl{A})\equiv \mathrm{Alg}(\mathrm{Lat}(\cl{A})).$$
 Whenever $\cl{A}=\mathrm{Ref}(\cl{A})$ we call $\cl{A}$ a \textbf{reflexive algebra}.

 A subspace $\cl{M}$ of $B(H_1,H_2)$ where $H_1, H_2$ are Hilbert spaces, is called a 
\textbf{ternary ring of operators (TRO)} if $$\cl{M}\cl{M}^*\cl{M}\subset \cl{M}.$$ 
In this case the spaces $[\cl{M}^*\cl{M}], [\cl{M}\cl{M^*}]$ are selfadjoint algebras. 
We call $\cl{M}$ \textbf{essential} if the algebras $[\cl{M}^*\cl{M}]^{-w^*}, 
[\cl{M}\cl{M^*}]^{-w^*}$ contain the identity operators. The following lemma is known. 
See for example \cite[8.5.32]{bm} or \cite[Lemma 3.1]{ele}.

\begin{lemma}\label{1.1} Let $\cl{C}, \cl{E}$ be von Neumann algebras 
acting on Hilbert spaces $H_1$ and $H_2$ respectively, $\theta : \cl{C}\rightarrow \cl{E}$ 
be a $*$-isomorphism and
$$\cl{M}=\{T\in B(H_1,H_2): TA=\theta (A)T \; \text{for\;\;all}\;\; A\in \cl{C}\}.$$ 
 The space $\cl{M}$ is an essential TRO.
\end{lemma}

\begin{definition}\label{1.1.d}\cite{ele}  Let $\cl{A}, \cl{B}$ be $w^*$ closed algebras 
acting on Hilbert spaces $H_1$ and $H_2$ respectively. If there exists a TRO 
$\cl{M}\subset B(H_1,H_2)$ such that $\cl{A}=[\cl{M}^*\cl{B}\cl{M}]^{-w^*}\;\; \text{and}\;\;
\cl{B}=[\cl{M}\cl{A}\cl{M}^*]^{-w^*}$ we write $\cl{A} 
\stackrel{\cl{M}}{\sim}\cl{B}.$ We say that the algebras $\cl{A}, \cl{B}$
 are \textbf{TRO equivalent} if there exists a TRO 
$\cl{M}$ such that $\cl{A} 
\stackrel{\cl{M}}{\sim} \cl{B}.$ 
\end{definition}

We also need the following main result of \cite{ele}.

\begin{theorem}\label{1.3} Two unital reflexive algebras $\cl{A}, \cl{B}$ are TRO equivalent 
if and only if there exists a $*$-isomorphism $\theta : \Delta(\cl{A})^\prime \rightarrow 
\Delta(\cl{B})^\prime $ such that $\theta (\mathrm{Lat}(\cl{A}))=\mathrm{Lat}(\cl{B}).$ If 
$\theta $ is as above and $\cl{M}=\{T: TA=\theta (A)T \; \text{for\;\;all}\;\; A\in 
\Delta(\cl{A})^\prime \}$ then $\cl{A} \stackrel{\cl{M}}{\sim}\cl{B}.$
\end{theorem}

We now define the category $_{\cl{A}}\mathfrak{M}$ for a unital dual operator
algebra $\cl{A}$ \cite{bm}. An object of  $_{\cl{A}}\mathfrak{M}$ is a Hilbert space $H$ 
for which there exists a unital algebraic homomorphism
$\alpha: \cl{A}\rightarrow B(H)$ which is  completely contractive and $w^*$-continuous.
We shall call such a map a \textbf{normal representation} of $\cl{A}.$ Throughout this work
we denote this object by $H$ or $(H,\alpha)$. If $(H_i,\alpha_i ), i=1,2$ are objects of 
$_{\cl{A}}\mathfrak{M}$ the space of homomorphisms $ \mathrm{ Hom}_{\cl{A}}(H_1,H_2)$ is the following: $$\mathrm{ Hom}_{\cl{A}}(H_1,H_2)=\{T\in B(H_1,H_2): T\alpha_1(A)= \alpha_2(A)T \text{\;for\;all\;}A\in \cl{A}\}.$$
Observe that the map $\alpha _i|_{\Delta (\cl{A})}$ is a $*$-homomorphism since
$\alpha _i$ is a contraction, \cite{bm}. We also define the category  $_{\cl{A}}\mathfrak{DM}$ 
which has the same objects as $_{\cl{A}}\mathfrak{M}$ but for every pair of objects 
$(H_i,\alpha_i ), i=1,2$ the space of homomorphisms 
$\mathrm{ Hom}_{\cl{A}}^{\mathfrak{D}}(H_1,H_2)$ is the following: $$\mathrm{ Hom}_{\cl{A}}^{\mathfrak{D}}(H_1,H_2)=\{T\in B(H_1,H_2): T\alpha_1(A)= \alpha_2(A)T \text{\;for\;all\;}A\in \Delta (\cl{A})\}.$$
 \;If the algebra $\cl{A}$ is a $W^*$-algebra the categories $_{\cl{A}}\mathfrak{M}$ and 
$_{\cl{A}}\mathfrak{DM}$ are the same. Also observe that $\mathrm{ Hom}_{\cl{A}}(H_1,H_2)\subset 
\mathrm{ Hom}_{\cl{A}}^{\mathfrak{D}}(H_1,H_2).$

\begin{definition}\label{1.2.d} \cite{ele2} Let $\cl{A}, \cl{B}$ be unital dual operator algebras 
and $\cl{F}:\;_{\cl{A}}\mathfrak{M} \rightarrow\; _{\cl{B}}\mathfrak{M}$ be a functor.
We say that the functor  $\cl{F}$ has a $\Delta$-extension if there exists a functor 
$\cl{G}: \;_{\cl{A}}\mathfrak{DM} \rightarrow \;_{\cl{B}}\mathfrak{DM} $ such that the following diagram is commutative:
\begin{align*}\begin{array}{clr} _{\cl{A}}\mathfrak{M}  & \hookrightarrow  
& _{\cl{A}}\mathfrak{DM} \\  \cl{F}\downarrow &  & \cl{G}\downarrow\;\;\;\;\;\;  \\  
_{\cl{B}}\mathfrak{M}  & \hookrightarrow  & _{\cl{B}}\mathfrak{DM}  \end{array} \end{align*}
\end{definition}

The following extends Rieffel's definition \cite{rif}: 

\begin{definition}\label{1.3.d}\cite{ele2}  Let $\cl{A}, \cl{B}$ be unital dual operator algebras and 
$\cl{F}: \;_{\cl{A}}\mathfrak{DM}\hookrightarrow \;_{\cl{B}}\mathfrak{DM}$  be a functor. 
We say that $\cl{F}$ is a $*$-\textbf{functor} if for every pair of objects  $H_1,H_2$ of 
$_{\cl{A}}\mathfrak{DM}$ every operator $F\in \mathrm{ Hom}_{\cl{A}}^{\mathfrak{D}}(H_1,H_2)$ satisfies $$\cl{F}(F^*)=\cl{F}(F)^*\in \mathrm{ Hom}_{\cl{B}}^{\mathfrak{D}}(\cl{F}(H_2),\cl{F}(H_1)).$$
\end{definition}

\begin{definition}\label{1.4.d}\cite{ele2}  Let $\cl{A}, \cl{B}$ be unital dual operator algebras. 
If there exists an equivalence functor $\cl{F}: \;_{\cl{A}}\mathfrak{M}  \rightarrow
 \;_{\cl{B}}\mathfrak{M}  $ which has a $\Delta$-extension as a $*$-functor 
implementing an equivalence between the categories $_{\cl{A}}\mathfrak{DM}  ,\;_{\cl{B}}
\mathfrak{DM},$ then $\cl{A}, \cl{B}$ are called \;\textbf{$\Delta $-equivalent} algebras.
\end{definition}

 The main theorem in \cite{ele2} which is a generalization of the main theorem of Rieffel 
\cite{rif}, is the following:

\begin{theorem}\label{1.6} Two unital dual operator algebras $\cl{A}, \cl{B}$ are 
$\Delta $-equivalent if and only if they have completely 
isometric normal representations $\alpha, \beta$ on Hilbert spaces such that the algebras  
$\alpha(\cl{A}), \beta(\cl{B})$ are TRO equivalent.
\end{theorem}

\section{Properties of the equivalence functors.}

In this section we fix unital dual operator algebras $\cl{A, B}$ and a functor $\cl{F}$
 implementing the equivalence of Theorem \ref{1.6}. We are going to  
 investigate properties of the functor $\cl{F},$ especially in case the algebras are 
reflexive.

 In \cite{ele2}, section 1, we proved that there exists 
$(H_0,\alpha _0)\in \;_{\cl{A}}\mathfrak{M}$ with corresponding object 
$(\cl{F}(H_0),\beta _0)\in \;_{\cl{B}}\mathfrak{M}$ such that $\alpha_0, \beta_0 $ 
are complete isometries and there exists an essential $w^*$-closed TRO $\cl{M}\subset B(H_0,\cl{F}(H_0))$ such that 
$\alpha _0(\cl{A}) 
\stackrel{\cl{M}}{\sim}\beta _0(\cl{B}).$ In what follows we identify $\cl{A}$ with 
$\alpha _0(\cl{A})$ and $\cl{B}$ with $\beta _0(\cl{B}).$ We denote by $\cl{U}$ 
and $\cl{V}$ the spaces $\cl{U}=[\cl{MA}]^{-w^*}, \cl{V}=[\cl{AM}^*]^{-w^*}$ which 
satisfy the following relations:
$$\cl{BUA}\subset \cl{U},\;\;\cl{AVB}\subset \cl{V},\;\;[\cl{VU}]^{-w^*}=\cl{A},\;\;[\cl{UV}]
^{-w^*}=\cl{B}.$$
Let us briefly recall from \cite{ele2}, section 2, the definition of a functor 
$\cl{F}_{\cl{U}}: \;_{\cl{A}}\mathfrak{M}\rightarrow 
\;_{\cl{B}}\mathfrak{M}$ and its $\Delta $-extension to  a $*$-functor 
$\cl{F}_{\cl{U}}: \;_{\cl{A}}\mathfrak{DM}\rightarrow 
\;_{\cl{B}}\mathfrak{DM}$ denoted by the same symbol.

For every  $(H,\alpha )\in \;_{\cl{A}}\mathfrak{M}$
the Hilbert space  $\cl{F}_{\cl U}(H)$ is the Hausdorff completion of the 
algebraic tensor product $\cl{U}\otimes H$ with respect to the following seminorm:
 $$\nor{\sum_{j=1}^m T_j\otimes x_j}_{\cl{F}_{\cl{U}}(H)}=
\sup_{S\in Ball(M_{n,1}(\cl{V})),n\in \mathbb{N}}\nor{\sum_{j=1}^m \alpha(ST_j)(x_j)}_{H^n}.$$
The representation $\cl F(\alpha)$ is defined by the formula
$$\beta=\cl F(\alpha) : 
\cl{B}\rightarrow B(\cl{F}_{\cl{U}}(H)): \beta (B)(T\otimes x)=BT\otimes x$$ 
where $B\in \cl{B}, T\in \cl{U}, x\in H$ and $T\otimes x$ is identified with its image in the quotient. 

Also for every
$H_1, H_2\in \;_{\cl{A}}\mathfrak{M}  $ we define a map
$$\cl{F}_{\cl{U}}: \mathrm{ Hom}_{\cl{A}}^\mathfrak{D}
(H_1,H_2)\rightarrow \mathrm{ Hom}_{\cl{B}}^\mathfrak{D}
(\cl{F}_{\cl{U}}(H_1),\cl{F}_{\cl{U}}(H_2))$$
by the formula
$$\cl{F}_{\cl{U}}(F)(M\otimes x)=M\otimes F(x)\;
\text{for\; all}\;F\in \mathrm{ Hom}_{\cl{A}}^\mathfrak{D}(H_1,H_2), M\in \cl{M}, x\in H_1.$$
The map $\cl{F}_{\cl{U}}(F)$ is well defined by this formula because 
$\cl{M}\otimes H_1$ is dense in $\cl{F}_{\cl{U}}(H_1),$ \cite[Corollary 2.5]{ele2}.

We will need the following theorem in \cite{ele2}:

\begin{theorem}\label{6.3} The functors $\cl{F}, \cl{F}_{\cl{U}}$ are equivalent as 
functors between the categories $_{\cl{A}}\mathfrak{M},\; _{\cl{B}}\mathfrak{M}
$ and their $\Delta$-extensions are equivalent as $*$-functors between the categories 
$_{\cl{A}}\mathfrak{DM}
, \;_{\cl{B}}\mathfrak{DM}.$
\end{theorem}

We now come to the concepts which will occupy us in this work.

\begin{definition}\label{6.1.d} Let $\cl{A}_1, \cl{B}_1$ be unital dual operator algebras.

(i) A functor $\cl{G}: \;_{\cl{A}_1}\mathfrak{M}
 \rightarrow \;_{\cl{B}_1}\mathfrak{M}$  is called \textbf{completely isometric (resp. normal)}
if for every pair of objects $H_1, H_2$ the map
$\cl{G}: \mathrm{ Hom}_{\cl{A}_1}(H_1,H_2) \to \mathrm{ Hom}_{\cl{B}_1}(\cl{G}(H_1),\cl{G}(H_2))$
is a complete isometry (resp. $w^*$-continuous).

(ii) We say that the functor $\cl{G}: \;_{\cl{A}_1}\mathfrak{M}
 \rightarrow \;_{\cl{B}_1}\mathfrak{M}
$ \textbf{respects isometries} if whenever $(H,\alpha )\in \;_{\cl{A}_1}\mathfrak{M}
$ is such that the map $\alpha : \cl{A}_1\rightarrow B(H)$ is a complete isometry
the corresponding map $\cl{G}(\alpha ): \cl{B}_1\rightarrow B(\cl{G}(H))$ is a complete isometry too.

(iii) We say that the functor $\cl{G}: \;_{\cl{A}_1}\mathfrak{M}
 \rightarrow \;_{\cl{B}_1}\mathfrak{M}
$ \textbf{respects reflexivity} if whenever $(H,\alpha )\in \;_{\cl{A}_1}\mathfrak{M}
$ is such that the map $\alpha : \cl{A}_1\rightarrow B(H)$ is a complete isometry and the algebra $\alpha (\cl{A}_1)$ is reflexive,  the map $\beta :\cl{B}_1\rightarrow B(\cl{G}(H))$ is a complete isometry and the algebra $\beta (\cl{B}_1)$ is reflexive, where $(\cl{G}(H),\beta )$ is the corresponding object of $_{\cl{B}_1}\mathfrak{M}
.$

(iv) A functor $\cl{G}: \;_{\cl{A}_1}\mathfrak{DM}
 \rightarrow \;_{\cl{B}_1}\mathfrak{DM}
$ is called a \textbf{lattice respecting} functor if for every object
$(H,\alpha )$ of $_{\cl{A}_1}\mathfrak{M}$
$$\cl{G}(\mathrm {Lat}(\alpha (\cl{A}_1)))=\mathrm {Lat}(\beta (\cl{B}_1))$$ where
$(\cl{G}(H),\beta )$ is the corresponding object of the category $_{\cl{B}_1}\mathfrak{M}.$
(Observe that $\mathrm {Lat}(\alpha (\cl{A}_1))\subset \mathrm{ Hom}_{\cl{A}_1}^\mathfrak{D}
(H,H)$ and $\mathrm {Lat}(\beta  (\cl{B}_1))\subset \mathrm{ Hom}_{\cl{B}_1}^\mathfrak{D}
(\cl{G}(H),\cl{G}(H)).$
\end{definition}

The following lemma is essentially due to Paschke; see for example \cite[8.5.23]{bm}.

\begin{lemma}\label{2.1}There exist partial isometries 
$\{W_k, k\in J\}\subset \cl{M}\;\; (\{V_k, k\in I\} \subset \cl{M})$ such that 
$W_k^*W_k\perp W_m^*W_m \;(V_kV_k^*\perp V_mV_m^*)$ for $k\neq m$ and 
$I_{H_0}=\sum_k\oplus W_k^*W_k \;(I_{\cl{F}(H_0)}=\sum_k\oplus V_kV_k^*).$
\end{lemma}

\begin{lemma}\label{6.4} The functor $\cl{F}_{\cl{U}}: \;_{\cl{A}}\mathfrak{M}
\rightarrow \;_{\cl{B}}\mathfrak{M}
$ respects isometries.
\end{lemma}
\textbf{Proof} Let $(H,\alpha )\in \;_{\cl{A}}\mathfrak{M}$ be such
that the representation $\alpha $ is a complete isometry. Suppose that $(\cl{F}_{\cl{U}}(H),\beta )$ is the corresponding object. We shall prove that the representation $\beta $ is a complete isometry too.

Let $n\in \bb N$. Fix a vector $y=(y_1,...,y_n)^t\in Ball(\cl{F}(H_0)^n)$ and a matrix
$(B_{ij})\in M_n(\cl{B})$. We recall from Lemma \ref{2.1} the partial
isometries $\{V_k, k\in I\} \subset \cl{M}.$ For $\epsilon >0$ there exists a subset $\{i_1,...,i_N\}\subset I$ such that
\begin{align*}
\|(B_{ij})(y)\|^2-\epsilon
=&\sum_{i=1}^n\nor{\sum_{k=1}^nB_{ik}(y_k)}^2-\epsilon \\
\leq &\sum_{i=1}^n\nor{\sum_{k=1}^nB_{ik}
\sum_{l=1}^NV_{i_l}V_{i_l}^*(y_k)}^2-\frac{\epsilon }{2}.
\end{align*}

Using again Lemma \ref{2.1} there exists a subset $\{j_1,...,j_m\}\subset I$ such that
\begin{align*}
&\sum_{i=1}^n\nor{\sum_{k=1}^nB_{ik}\sum_{l=1}^NV_{i_l}V_{i_l}^*(y_k)}^2-\frac{\epsilon }{2}\\
\leq & \sum_{i=1}^n\nor{\sum_{t=1}^mV_{j_t}V_{j_t}^*\left(\sum_{k=1}^nB_{ik}
\sum_{l=1}^NV_{i_l}V_{i_l}^*(y_k)\right)}^2-\frac{\epsilon }{4}.
\end{align*}
Since the projections $(V_{j_t}V_{j_t}^*)$ are mutually orthogonal it follows that
\begin{align*}
&\|(B_{ij})(y)\|^2-\epsilon \leq \sum_{i=1}^n\sum_{t=1}^m \nor{V_{j_t}^*\left(\sum_{k=1}^nB_{ik}\sum_{l=1}^NV_{i_l}V_{i_l}^*(y_k)\right)}^2-\frac{\epsilon }{4}\\
=&\|(V^* \oplus ...\oplus V^* )(B_{ij}) (U \oplus ...\oplus U)(z)\|^2-\frac{\epsilon }{4},
\end{align*} where
$$z=(V_{i_1}^* (y_1),...,V_{i_N}^*(y_1),V_{i_1}^*(y_2),...,V_{i_N}^*(y_2),...,V_{i_N}^*(y_n))^t,$$
$$V=(V_{j_1},...,V_{j_m}), U=(V_{i_1},...,V_{i_N}).$$

Observe that
\begin{align*}
&\|z\|^2=\sum_{l=1}^n \sum_{k=1}^N \|V_{i_k}^*(y_l)\|^2
=\sum_{l=1}^n\sum_{k=1}^N \|V_{i_k}V_{i_k}^*(y_l)\|^2\\
=&\sum_{l=1}^n\nor{\sum_{k=1}^N V_{i_k}V_{i_k}^*(y_l)}^2
\leq \sum_{l=1}^n \|y_l\|^2=\|y\|^2 \leq 1.
\end{align*}

 It follows that
\begin{align*}
\|(B_{ij})(y)\|^2-\epsilon& \leq \|(V^*\oplus ...\oplus V^*)(B_{ij}) (U\oplus ...\oplus U)\|^2-\frac{\epsilon }{4}\\
=& \|\alpha \left((V^*\oplus ...\oplus V^*)(B_{ij}) (U\oplus ...\oplus U)\right)\|^2-\frac{\epsilon }{4}
\end{align*}
(the last equality holds because the map $\alpha $ is a complete isometry).
We can find vectors $x_{lk}\in H$ such that the vector $$x=(x_{11},...,x_{N1},x_{12},...,x_{N2},...,x_{Nn})^t\in H^{Nn}$$
has norm one and \begin{align*}
\|(B_{ij})(y)\|^2-\epsilon \leq \|\alpha ((V^*\oplus ...\oplus V^*)(B_{ij}) (U\oplus ...\oplus U))(x)\|^2.
\end{align*}
Thus
\begin{equation}\label{eq62}
\|(B_{ij})(y)\|^2-\epsilon\leq \sum_{s=1}^n\sum_{r=1}^m\nor{\sum_{k=1}^n\sum_{l=1}^N \alpha (V_{j_r}^*B_{sk}V_{i_l})(x_{lk})}^2.
 \end{equation}

Let
$$\omega _k=\sum_{l=1}^NV_{i_l}\otimes x_{lk}\in \cl{F}_{\cl{U}}(H),\; k=1,...,n.$$
We have
\begin{align*}&\sum_{k=1}^n \|\omega _k \|^2_{\cl{F}_{\cl{U}}(H)}=
\sum_{k=1}^n \sup_{S \in Ball(M_{p,1}(\cl{V})), p\in \mathbb{N}}
\nor{ \sum_{l=1}^N\alpha (SV_{i_l})(x_{lk}) }^2 \\
=&\sum_{k=1}^n \sup_{S\in Ball(M_{p,1}(\cl{V})), p \in \mathbb{N}}
\nor{\alpha (S(V_{i_1}...V_{i_N}))\left( \begin{array}{clr} x_{1k}\\
\vdots \\
x_{Nk} \end{array}\right)}^2 \leq \sum_{k=1}^n\nor{\left(\begin{array}{clr}x_{1k}\\
\vdots\\
x_{Nk}\end{array}\right)}^2\leq 1.
\end{align*}
The inequality is a consequence of the fact that $\|(V_{i_1}...V_{i_N})\|=1$ and the map
$\alpha $ is a complete isometry.

It follows that
\begin{align*}  &\|\beta ((B_{ij}))\|^2 \geq \|\beta ((B_{ij}))(\omega_1...\omega _n)^t\|^2\\=&\sum_{s=1}^n\nor{\sum_{k=1}^n\beta (B_{sk})(\omega _k)}^2_{\cl{F}_{\cl{U}}(H)}= \sum_{s=1}^n\nor{\sum_{k=1}^n\sum_{l=1}^NB_{sk}V_{i_l}\otimes x_{lk}}^2_{\cl{F}_{\cl{U}}(H)}\\=&\sum_{s=1}^n\sup_{S\in Ball(M_{t,1}(\cl{V})), t\in \mathbb{N}}\nor{\sum_{k=1}^n\sum_{l=1}^N \alpha (SB_{sk}V_{i_l})(x_{lk})}^2
\end{align*}

Observe that $$S=(V_{j_1}^*...V_{j_m}^*)^t\in Ball(M_{m,1}(\cl{V}))$$ so we have 
\begin{align*} \|\beta ((B_{ij}))\|^2\geq & \sum_{s=1}^n\nor{\sum_{k=1}^n\sum_{l=1}^N
\alpha ((V_{j_1}^*...V_{j_m}^*)^tB_{sk}V_{i_l})(x_{lk})}^2\\=&\sum_{s=1}^n\sum_{r=1}^
m\nor{\sum_{k=1}^n\sum_{l=1}^N \alpha (V_{j_r}^*B_{sk}V_{i_l})(x_{lk})}^2. \end{align*}

From inequality (\ref{eq62}) it follows that $\|(B_{ij})(y)\|^2-\epsilon \leq \|\beta 
((B_{ij}))\|^2,$ hence $\|(B_{ij})\| \leq \|\beta ((B_{ij}))\|.$ Since $\beta$ is a 
complete contraction we have equality: $\|(B_{ij})\| = \|\beta ((B_{ij}))\|.$ $\qquad  \Box$

\bigskip

Combining this lemma and Theorem \ref{6.3} we obtain the next theorem.

\begin{theorem}\label{6.5} Every functor implementing the equivalence of Theorem \ref{1.6}
 respects isometries.
\end{theorem}

\begin{lemma}\label{6.6} Let $(H,\alpha )\in \;_{\cl{A}}\mathfrak{M}
$ and $(\cl{F}_{\cl{U}}(H), \beta )$ be the corresponding object. Then 
$$\cl{F}_{\cl{U}}(\mathrm {Lat}(\alpha (\cl{A})))\subset 
\mathrm {Lat}(\beta  (\cl{B})). $$ 
\end{lemma}
\textbf{Proof}
  
Suppose that $L$ is a projection of the lattice of $\alpha (\cl{A}).$ We shall prove that
$\cl{F}_{\cl{U}}(L)\in \mathrm {Lat}(\beta (\cl{B})).$
The operator $\cl{F}_{\cl{U}}(L)$ is a projection
because $\cl{F}_{\cl{U}}$ is a $*$-functor (\cite[Theorem 2.10]{ele2}). If $B\in \cl{B}, M\in \cl{M}$ and $x\in H$ then
$$\beta (B)\cl{F}_{\cl{U}}(L)(M\otimes x)=\beta (B)(M\otimes L(x))=BM\otimes L(x).$$
By \cite[Lemma 2.4]{ele2} we have
$$BM\otimes L(x)\in [N\otimes L(y): N\in \cl{M}, y\in H]^{-}
=[\cl{F}_{\cl{U}}(L)(N\otimes y): N\in \cl{M}, y\in H]^{-}.$$
It follows that
$$\beta (B)\cl{F}_{\cl{U}}(L)(M\otimes x)\in \cl{F}_{\cl{U}}(L)(\cl{F}_{\cl{U}}(H)).$$
Since the space $\cl{M}\otimes H$ is dense in $\cl{F}_{\cl{U}}(H),$ 
\cite[Corollary 2.5]{ele2}, we obtain
$$\beta (B)\cl{F}_{\cl{U}}(L)(z)\in \cl{F}_{\cl{U}}(L)(\cl{F}_{\cl{U}}(H))\;
\text{for\;all\;}z\in \cl{F}_{\cl{U}}(H), B\in \cl{B}.$$
This shows that $\cl{F}_{\cl U}(L)\in \mathrm{Lat}(\beta(\cl B)).$ 
We proved that $$ \cl{F}_{\cl{U}}(\mathrm {Lat}(\alpha (\cl{A})))\subset 
\mathrm {Lat}(\beta  (\cl{B})).  \qquad  \Box$$

\begin{theorem}\label{6.7} Every functor implementing the equivalence of Theorem \ref{1.6} is a lattice respecting functor.
\end{theorem}
\textbf{Proof} Let $(H,\alpha )\in \;_{\cl{A}}\mathfrak{M}
$ and $(\cl{F}(H),\cl{F}(\alpha ))$ $\in \;_{\cl{B}}\mathfrak{M}
$ be the corresponding object. Since the functor $\cl{F}$ is equivalent to $\cl{F}_{\cl{U}}$ 
it follows from the last lemma that  $$\cl{F}(\mathrm {Lat}(\alpha (\cl{A})))\subset 
\mathrm {Lat}(\cl{F}(\alpha )(\cl{B})).$$
Suppose that $\cl{G}$ is the inverse functor of $\cl{F}$ which maps 
 $(\cl{F}(H),\cl{F}(\alpha ))$ $\in \;_{\cl{B}}\mathfrak{M}$ to 
$(\cl{G}\cl{F}(H),\cl{G}\cl{F}(\alpha ))$ $\in \;_{\cl{A}}\mathfrak{M}.$
 By the same argument $$\cl{G}(\mathrm {Lat}(\cl{F}(\alpha) (\cl{B})))\subset 
\mathrm {Lat}(\cl{G}\cl{F}(\alpha )(\cl{A})).$$
If $\cl{F}( \mathrm {Lat}(\alpha (\cl{A} )))$ is strictly contained in 
$\mathrm {Lat}(\cl{F}(\alpha )(\cl{B}))$ then 
$\cl{G}\cl{F}(\mathrm {Lat}(\alpha (\cl{A})))$ is strictly contained in 
$\mathrm {Lat}(\cl{G}\cl{F}(\alpha )(\cl{A})).$ The functor $\cl{GF}$ is equivalent 
to the identity functor of the category $\;_{\cl{A}}\mathfrak{M}.$ So there exists 
unitary $U\in \mathrm{ Hom}_{\cl{A}}(\cl{GF}(H),H)$ satisfying $U^*FU=\cl{GF}(F)$
 for all $F\in \mathrm{ Hom}_{\cl{A}}^{\mathfrak{D}}(H,H).$ So 
$$\cl{GF}(\mathrm {Lat}(\alpha (\cl{A}))=U^* \mathrm {Lat}(\alpha (\cl{A}))U =
 \mathrm {Lat}(\cl{GF}(\alpha )(\cl{A})).$$ This is a contradiction, and hence we have 
the equality $$\cl{F}(\mathrm {Lat}(\alpha (\cl{A})))= 
\mathrm {Lat}(\cl{F}(\alpha )(\cl{B})). \qquad \Box$$

\begin{theorem}\label{6.11.a} Let $(H,\alpha )\in \;_{\cl{A}}\mathfrak{M}$ be such
that $\alpha $ is a complete isometry.
If $(\cl{F}(H),\beta )\in \;_{\cl{B}}\mathfrak{M}$ is the
corresponding object then $\beta $ is a complete isometry and the
algebras $\alpha(\cl{A}), \beta(\cl{B}) $ are TRO equivalent.
\end{theorem}
\textbf{Proof} By Theorem \ref{6.5}, $\beta $ is a complete
isometry. We denote by $\sigma $ the map $$\cl{F}:
\mathrm{ Hom}_{\cl{A}}^\mathfrak{D}  (H,H)=\alpha (\Delta
(\cl{A}))^\prime\rightarrow \beta (\Delta
(\cl{B}))^\prime=\mathrm{ Hom}_{\cl{B}}^\mathfrak{D}
(\cl{F}(H),\cl{F}(H))$$ which is a $*$-isomorphism. By Lemma
\ref{1.1} the space
$$\cl{Y}=\{N\in B(H,\cl{F}(H)): NA=\sigma (A)N\;\text{for\;all\;}A\in
\alpha (\Delta (\cl{A}))^\prime\}$$
is an essential TRO. In the sequel, if $K$ is a Hilbert space,
$T$ is an operator on $K$ and $\cl C\subset B(K)$  we
denote by $K^\infty $ the countably infinite direct sum
$K\oplus K \oplus \ldots$, by $T^\infty\in B(K^\infty)$ the
operator $T\oplus T \oplus ...,$ and by $\cl{C}^\infty $ the
set $\{C^\infty : C\in \cl{C}\}.$

The map $\alpha ^\infty : \cl{A}\rightarrow B(H^\infty )$ given by
$\alpha ^\infty (A)=\alpha (A)^\infty $ is a normal representation
so $(H^\infty, \alpha ^\infty) \in \;_{\cl{A}}\mathfrak{M}$ hence
$(\cl{F}(H^\infty), \cl{F}(\alpha ^\infty)) \in
\;_{\cl{B}}\mathfrak{M}.$ Let $U\in \mathrm{ Hom}_{\cl{B}}(\cl{F}(H^\infty
), \cl{F}(H)^\infty )$ be a unitary (see \cite[Lemma 3.2]{ele2}). This defines a unitary equivalence between the algebras
$$\mathrm{ Hom}_{\cl{B}}^{\mathfrak{D}}(\cl{F}(H^\infty ), \cl{F}(H^\infty)
)=(\cl{F}(\alpha ^\infty )(\Delta (\cl{B})))^\prime$$ and $$\mathrm{ Hom}_{\cl{B}}^
{\mathfrak{D}}(\cl{F}(H)^\infty ,
\cl{F}(H)^\infty)=(\beta (\Delta (\cl{B}))^\infty)^\prime.$$ This equivalence
 maps the invariant projection lattice of the
algebra $\cl{F}(\alpha ^\infty )(\cl{B})$ onto the lattice of
$\beta (\cl{B})^\infty.$ The functor $\cl{F}$ defines  a
$*$-isomorphism between the spaces
$$\mathrm{ Hom}_{\cl{A}}^{\mathfrak{D}}(H^\infty ,H^\infty)=(\alpha (\Delta
(\cl{A}))^\infty)^\prime\;\;\text{and}\;\; \mathrm{ Hom}_{\cl{B}}^{\mathfrak{D}}(\cl{F}
(H^\infty ), \cl{F}(H^\infty)
)$$ which by Theorem \ref{6.7} maps the lattice of the algebra
$\alpha (\cl{A})^\infty $ onto the lattice of the algebra
$\cl{F}(\alpha ^\infty )(\cl{B}).$ Composing with the unitary $U$ 
we obtain a $*$-isomorphism $$\theta :
(\alpha (\Delta (\cl{A}))^\infty)^\prime\rightarrow (\beta (\Delta
(\cl{B}))^\infty)^\prime $$ such that $\theta (\mathrm{Lat}(\alpha
(\cl{A})^\infty )) = \mathrm{Lat}(\beta  (\cl{B})^\infty) $ and
which satisfies $$\theta ((F_{ij})_{i,j})=(\sigma
(F_{ij}))_{i,j}\;\text{for \; all \;}(F_{ij})_{i,j}\in (\alpha
(\Delta (\cl{A}))^\infty)^ \prime.$$ From this we conclude that the space 
$$\cl{X}=\{(T_{ij})\in B(H^\infty, \cl{F}(H)^\infty ): (T_{ij})(F_{ij})
=\theta ((F_{ij}))(T_{ij})\;\forall \;(F_{ij}) \in (\alpha (\Delta
(\cl{A}))^\infty)^ \prime\},$$
equals $\cl Y^\infty$. Since the algebras
$\alpha(\cl{A})^\infty , \beta(\cl{B})^\infty  $ are reflexive
(see for example \cite[A.1.5]{bm}) by Theorem \ref{1.3} we have
$$ \alpha (\cl{A})^\infty  \stackrel{\cl{X}}{\sim} \beta (\cl{B})^\infty \Rightarrow 
 \alpha (\cl{A})  \stackrel{\cl{Y}}{\sim} \beta (\cl{B}).\qquad \Box $$

\begin{theorem}\label{6.11} Every functor implementing the equivalence of
Theorem \ref{1.6} respects reflexivity.
\end{theorem}

\textbf{Proof} Let $(H,\alpha )\in \;_{\cl{A}}\mathfrak{M}  $ be such that
$\alpha $ is a complete isometry and the algebra $\alpha (\cl{A})$ is reflexive.
Suppose that $(\cl{F}(H),\beta )\in \;_{\cl{B}}\mathfrak{M}  $ is the corresponding
object. By Theorem \ref{6.5} $\beta $ is a complete isometry. By the above theorem
the algebras $\alpha(\cl{A}), \beta(\cl{B}) $ are TRO equivalent. Since
$\alpha (\cl{A})$ is reflexive so is $\beta  (\cl{B})$ (\cite[Remark 2.7]
{ele}). $\qquad  \Box$

\section{Applications and examples.}

We present some definitions and concepts used in this section.
A \textbf{commutative subspace lattice (CSL)} is a projection lattice
$\cl{L}$ whose elements commute;
the algebra $\mathrm{Alg}(\cl{L})$ is called a \textbf{CSL algebra}.
In the special case where $\cl{L}$ is totally ordered we call
 $\cl{L}$ a \textbf{nest} and the algebra $\mathrm{Alg}(\cl{L})$
a \textbf{nest algebra}. CSL algebras are of course reflexive.
When $\cl{A}$ is a CSL algebra there exists a smallest $w^*$-closed algebra
contained in $\cl{A},$ which contains the diagonal $\Delta (\cl{A})$ and whose
reflexive hull is $\cl{A}$ \cite{arv}, \cite{st}. We denote this algebra by $\cl{A}_{min}.$ Whenever $\cl{A}= \cl{A}_{min}$ we call $\cl{A}$ \textbf{synthetic}. The first example of a nonsynthetic CSL algebra was given in \cite{arv}.

\begin{proposition}\label{7.1} If $\cl{A}$ is a CSL algebra which is $\Delta $-equivalent to a
unital dual operator algebra $\cl{B}$ then there exists a
completely isometric normal representation $\beta$ of $\cl{B}$
such that the algebras $\cl{A}$ and $\beta (\cl{B})$ are TRO
equivalent. It follows that the algebra $\beta (\cl{B})$ is a CSL
algebra too.
\end{proposition}
\textbf{Proof} Suppose that $\cl{F}: \;_{\cl{A}}\mathfrak{M}
\rightarrow \;_{\cl{B}}\mathfrak{M}$ is an equivalence  functor
which has a $\Delta$-extension to an equivalence
$*$-functor between the categories $_{\cl{A}}\mathfrak{DM}$ and
$_{\cl{B}}\mathfrak{DM}$. Also suppose that $\cl{A}\subset B(H)$
and let $(\cl{F}(H),\beta )$ be the object corresponding to identity representation of $\cl{A}.$ 
 By Theorem \ref{6.11} $\beta $ is a complete isometry and the algebra 
$\beta (\cl{B})$ is reflexive. Also by Theorem \ref{6.11.a} the algebras $\cl{A}, \beta
(\cl{B})$ are TRO equivalent. Now Theorem \ref{1.3} shows that
the lattice $\mathrm{Lat}(\beta (\cl{B}))$ is a CSL. $\qquad \Box$

\medskip 

Although $\Delta $-equivalent algebras are not necessarily TRO-equivalent, even 
when they are reflexive, in CSL algebras this is indeed the case:

\begin{theorem}\label{7.2} Two CSL algebras are $\Delta $-equivalent if and only if they are TRO-equivalent.
\end{theorem}
\textbf{Proof} TRO-equivalent algebras are $\Delta$-equivalent
(Theorem \ref{1.6}). For the converse, suppose that the CSL
algebras $\cl{A}$ and $\cl{B}$ are $\Delta$-equivalent. By the
previous proposition there exists a completely isometric normal
representation $\beta $ of $\cl{B}$ such that the algebras
$\cl{A}$ and $\beta (\cl{B})$ are TRO equivalent. Since $\beta(\cl B)$
is a CSL algebra, as just shown, it is easily
checked that $\beta(\mathrm {Lat}(\cl{B}))=\mathrm
{Lat}(\beta(\cl{B}))$ and $\beta (\Delta (\cl{B})')=(\Delta (\beta
(\cl{B})))'$. It follows (Theorem \ref{1.3}) that the
algebras $\cl{B}$ and $\beta (\cl{B})$ are TRO equivalent. The
conclusion is a consequence of the fact that TRO equivalence is an
equivalence relation (Theorem 2.3 in \cite{ele}). $\qquad \Box$

\begin{remark}\label{7.3}
\em{(i) Suppose that $\cl{A}$ and $\cl{B}$ are separably acting CSL
algebras with continuous or totally atomic lattices. By the
previous theorem and \cite[Theorem 5.7]{ele} they are
$\Delta$-equivalent if and only if they have isomorphic lattices. 
In general, two separably acting CSL algebras are $\Delta$-equivalent
if and only if their lattices are isomorphic through a lattice isomorphism
which ``respects continuity" \cite[Theorem 5.7]{ele}.

(ii) If two nests are isomorphic, their nest algebras are not always
$\Delta$-equivalent, even if they have isomorphic diagonals (see example \ref{7.7}).}
\end{remark}

\begin{proposition}\label{7.4}
If $\cl{A}$ is a nonsynthetic CSL algebra there exists no isometric normal
representation $\alpha : \cl{A}_{min}\to B(H)$ such that $\alpha
(\cl{A}_{min})$ is a CSL algebra. It follows from Proposition
\ref{7.1} that the algebra $\cl{A}_{min}$ cannot be
$\Delta$-equivalent to any CSL algebra.
\end{proposition}
\textbf{Proof} Let $\cl{A}$ be a nonsynthetic CSL
algebra, $H$  a Hilbert space and $\alpha :
\cl{A}_{min}\rightarrow B(H)$ be a $w^*$-continuous isometric
homomorphism such that $\cl{B}\equiv \alpha (\cl{A}_{min})$ is a
CSL algebra. Since $\cl{A}$ equals $ \mathrm {Ref}(\cl{A}_{min})$
and $\mathrm {Lat}(\cl{A})\subset \cl{A}_{min}$ we can check that
$\alpha (\mathrm {Lat}(\cl{A}))=\mathrm {Lat}(\cl{B}).$
From \cite[Theorem 4.7]{ele} the algebra $\cl{B}$ is not synthetic,
so the algebra $\cl{B}_{min}$ is strictly contained in $\cl{B}$.
Thus the algebra $\alpha ^{-1}(\cl{B}_{min})$ is
strictly contained in $\cl{A}_{min}.$ This is a contradiction
because $\Delta (\cl{A})\subset \alpha ^{-1}(\cl{B}_{min}),
\mathrm {Lat}(\cl{A})=\mathrm{Lat}(\alpha ^{-1}(\cl{B}_{min}))$
and the algebra $\cl{A}_{min}$ is the smallest $w^*$-closed
subalgebra of $\cl{A}$ with these properties \cite{arv}, \cite{st}. $\qquad
\Box$

\bigskip

 We will prove that similar nest algebras \cite{dav} have
equivalent categories. So we fix nests $\cl{N}_1, \cl{N}_2$
acting on the separable Hilbert spaces $H_1, H_2$ respectively and
an order
isomorphism $\theta : \cl{N}_1\rightarrow \cl{N}_2$ which preserves
dimension of intervals. We say that an invertible operator $T\in B(H_1,H_2)$
{\em implements} $\theta $ if $\theta (N)$ is the projection onto
$TN(H_1)$ for all $N\in \cl{N}_1.$ Define the spaces
$$\cl{U}=\{T\in B(H_1,H_2): \theta (N)^\bot TN=0\;\text{for\;all\;}N\in \cl{N}_1\},$$
 $$\cl{V}=\{S\in B(H_2,H_1): N^\bot S\theta (N)=0\;\text{for\;all\;}N\in \cl{N}_1\}.$$
 If $\cl{A}=\mathrm {Alg}(\cl{N}_1)$ and
$\cl{B}=\mathrm{Alg}(\cl{N}_2)$, one verify easily that 
$$\cl{V}\cl{U}\subset \cl{A},\; \cl{U}\cl{V}\subset \cl{B},\; \cl{B}\cl{U}\cl{A}\subset \cl{U},\; 
\cl{A}\cl{V}\cl{B}\subset \cl{V}.$$

\medskip

We will need the Similarity Theorem \cite[Theorem 13.20]{dav}.

\begin{theorem}[Davidson]\label{7.5}
For every $\epsilon >0$ there exists an invertible operator $T$
which implements $\theta $ such that $\|T\|<1+\epsilon , \|T^{-1}\|<1+\epsilon .$
\end{theorem}

\begin{proposition}\label{7.6} There exists an equivalence functor 
$ \cl{F}_{\cl{U}} : \;_{\cl{A}}\mathfrak{M} \rightarrow \; _{\cl{B}}\mathfrak{M}  $ 
which is normal and completely isometric.
\end{proposition}
\textbf{Proof} Let $(H,\alpha )\in \;_{\cl{A}}\mathfrak{M}.$ We
define the following seminorm on the algebraic tensor product
$\cl{U}\otimes H:$
$$\nor{\sum_{i=1}^nT_i\otimes x_i}=\sup_{S\in Ball(\cl{V})}\nor{\sum_{i=1}^n\alpha (ST_i)(x_i)}.$$
This seminorm satisfies the parallelogram identity. Let
$\cl{F}_{\cl{U}}(H)$ be the corresponding Hausdorff completion of
$\cl U\otimes H$ and identify every $T\otimes x$ with its image in
$\cl{F}_{\cl{U}}(H)$. If $B\in \cl{B},
T_1,...,T_m\in \cl{U}, x_1,...,x_m\in H$ we can check that
$$\nor{\sum_{j=1}^mBT_j\otimes x_j}_{\cl{F}_{\cl{U}}(H)}\leq \|B\|
\nor{\sum_{j=1}^mT_j\otimes x_j}_{\cl{F}_{\cl{U}}(H)}.$$
The map $T\otimes x\rightarrow BT\otimes x, B\in \cl{B}, T\in \cl{U}, x\in H$ extends to a map 
$\beta (B)\in B(\cl{F}_{\cl{U}}(H))$ and clearly $\beta $ is a unital algebraic homomorphism.
 As in \cite[Proposition 2.7]{ele2} we can prove that the map 
$\beta : \cl{B}\rightarrow B(\cl{F}_{\cl{U}}(H))$ is $w^*$-continuous. It 
follows from \cite[Corollary 20.17]{dav} that $\beta $ is a complete contraction and hence 
$(\cl{F}_{\cl{U}}(H),\beta )\in \;_{\cl{B}}\mathfrak{M}.$ In order to define a functor $\cl{F}_{\cl{U}}: \;_{\cl{A}}\mathfrak{M}  \rightarrow  \;_{\cl{B}}\mathfrak{M}  $ we have to define the map $$\cl{F}_{\cl{U}}: \mathrm{ Hom}_{\cl{A}}(H_1,H_2)\rightarrow  \mathrm{ Hom}_{\cl{B}}(\cl{F}_{\cl{U}}(H_1),\cl{F}_{\cl{U}}(H_2))$$ for every pair of objects $(H_j,\alpha _j), j=1,2.$ If $F\in \mathrm{ Hom}_{\cl{A}}(H_1,H_2),$ then \begin{align*}&\nor{\sum_iT_i\otimes F(x_i)}_{\cl{F}_{\cl{U}}(H_2)}=\sup_{S\in Ball(\cl{V})}\nor{\sum_i\alpha _2(ST_i)F(x_i)}\\=&\sup_{S\in Ball(\cl{V})}\nor{F\sum_i\alpha _1(ST_i)(x_i)}\leq \|F\| \nor{\sum_iT_i\otimes x_i}_{\cl{F}_{\cl{U}}(H_1)}.
\end{align*}
 So we can define a map $\cl{F}_{\cl{U}}(F)\in B(\cl{F}_{\cl{U}}(H_1),\cl{F}_{\cl{U}}(H_2))$ by the formula  $$\cl{F}_{\cl{U}}(F)(T\otimes x)=T\otimes F(x), T\in \cl{U}, x\in H_1.$$
 It is easy to check that $\cl{F}_{\cl{U}}(F)\in \mathrm{ Hom}_{\cl{B}}(\cl{F}_{\cl{U}}(H_1),\cl{F}_{\cl{U}}(H_2)).$ The definition of the functor $\cl{F}_{\cl{U}}$ is complete. Symmetrically \;we can define a \;\;functor  $\cl{F}_{\cl{V}}:\;\;  _{\cl{B}}\mathfrak{M}  \rightarrow  \;_{\cl{A}}\mathfrak{M}  .$

 Now fix $(H,\alpha )\in \;_{\cl{A}}\mathfrak{M}  $ with corresponding object $(\cl{F}_{\cl{U}},\beta ).$ If $S_i\in \cl{V}, T_i\in \cl{U}, x_i\in H, i=1,...,r,$ then

\begin{align*}&\nor{\sum_{i=1}^rS_i\otimes (T_i\otimes x_i)}_{\cl{F}_
{\cl{V}}\cl{F}_{\cl{U}}(H)}=\sup_{U\in Ball(\cl{U})}\nor{\sum_{i=1}^r\beta 
(US_i)(T_i\otimes x_i)}_{\cl{F}_{\cl{U}}(H)}\\=&\sup_{U\in Ball(\cl{U})}\nor
{\sum_{i=1}^rUS_iT_i\otimes x_i}_{\cl{F}_{\cl{U}}(H)}=\sup_{U\in Ball(\cl{U})}
\sup_{V\in Ball(\cl{V})}\nor{\sum_{i=1}^r\alpha (VU)\alpha (S_iT_i)(x_i)}_H\\\leq & 
\nor{\sum_{i=1}^r\alpha (S_iT_i)(x_i)}_H
\end{align*}
 By Theorem \ref{7.5} for arbitrary $\epsilon >0$ we can choose $T\in \cl{U}$ such that $T^{-1}\in \cl{V}$ and $\|T\|<1+\epsilon , \|T^{-1}\|<1+\epsilon .$ By the definition of the norm $\|\cdot\|_{\cl{F}_{\cl{V}}\cl{F}_{\cl{U}}(H)}$ we have
\begin{align*}&\nor{\sum_{i=1}^rS_i\otimes (T_i\otimes x_i)}_{\cl{F}_{\cl{V}}\cl{F}_{\cl{U}}(H)}\geq \nor{\sum_{i=1}^r\alpha \left(\frac{T^{-1}}{\|T^{-1}\|}\frac{T}{\|T\|}S_iT_i\right)(x_i)}\\\geq &\frac{1}{(1+\epsilon )^2}\nor{\sum_{i=1}^r\alpha (S_iT_i)(x_i)}
\end{align*}
 Letting $\epsilon \rightarrow 0$ we obtain $$\nor{\sum_{i=1}^rS_i\otimes (T_i\otimes x_i)}_{\cl{F}_{\cl{V}}\cl{F}_{\cl{U}}(H)}\geq \nor{\sum_{i=1}^r\alpha (S_iT_i)(x_i)}$$ and hence equality holds. It follows that we can define a unitary $$U_H: \cl{F}_{\cl{V}}\cl{F}_{\cl{U}}(H)\rightarrow H: U_H(S\otimes (T\otimes x))=\alpha (ST)(x), S\in \cl{V}, T\in \cl{U}, x\in H.$$
 One can now easily check that the family of unitaries $\{U_H: H\in \;_{\cl{A}}\mathfrak{M}  \}$ 
implements the required equivalence. The proofs of the facts that the functor 
$ \cl{F}_{\cl{U}}$ is  normal and completely isometric are similar to \cite[Lemmas 3.4, 
3.5]{ele2} so we omit them.$\qquad \Box$

\medskip

 In spite of the last proposition, we show in the following example that the similarity of nest algebras does not imply $\Delta $-equivalence even in the case of isomorphic diagonals.

\begin{example}\label{7.7}\em{ In \cite[Examples 13.25, 13.22]{dav} 
there exist similar nests $\cl{N}_1, \cl{N}_2$ acting 
on separable Hilbert spaces $H_1,H_2$ respectively, such that the 
algebra $\cl{N}_1^{\prime\prime}$ is a totally atomic maximal abelian 
selfadjoint algebra (masa) and the algebra $\cl{N}_2^{\prime\prime}$ is 
a masa with the property that the algebra $\cl{N}_2^{\prime\prime}|_{N(H_2)}$ 
has nontrivial continuous part for every nonzero projection $N\in \cl{N}_2.$ 
We define the nests $$\cl{M}_1=\{0\oplus N: N\in \cl{N}_1\}\cup \{N\oplus H_1: 
N\in \cl{N}_2\}\subset B(H_2\oplus H_1)$$ $$\cl{M}_2=\{0\oplus N: N\in 
\cl{N}_2\}\cup \{N\oplus H_2: N\in \cl{N}_1\}\subset B(H_1\oplus H_2).$$ 
These nests are similar too.

 So by the previous proposition if $\cl{A}=\mathrm {Alg}(\cl{M}_1), 
\cl{B}=\mathrm {Alg}(\cl{M}_2)$ the categories  $_{\cl{A}}\mathfrak{M}  
,\; _{\cl{B}}\mathfrak{M}  $ are equivalent. Observe that the diagonals 
of these algebras are isomorphic because $\Delta (\cl{A})=\cl{N}_2^{\prime\prime}
\oplus \cl{N}_1^{\prime\prime}$ and $\Delta (\cl{B})=\cl{N}_1^{\prime\prime}
\oplus \cl{N}_2^{\prime\prime}$. The algebras $\cl{A}$ and $\cl{B}$ are not 
$\Delta $-equivalent because if they were by Theorem  \ref{7.2} they would be 
TRO equivalent. So by Theorem \ref{1.3} there would exist a $*$-isomorphism 
$$\pi :\Delta(\cl{A})\rightarrow  \Delta(\cl{B})\;\text{such\;that\;}\pi 
(\cl{M}_1)=\cl{M}_2.$$ Since the diagonals are masas this map is unitarily 
implemented \cite[Theorem 9.3.1]{kr}. Now there are two possibilities: either 
$\pi (0\oplus I_{H_1})=0\oplus N$ for some $N\in \cl{N}_2$ in which case the 
algebras $\cl{N}_1^{\prime\prime}$ and $\cl{N}_2^{\prime\prime}|_{N(H_2)}$ 
are unitarily equivalent, or $\pi (0\oplus I_{H_1})=M\oplus I_{H_2}$ for 
some $M\in \cl{N}_1$ in which case the algebras $\cl{N}_1^{\prime\prime}$ 
and $\cl{N}_1^{\prime\prime}|_{M(H_1)}\oplus \cl{N}_2^{\prime\prime}$ are 
unitarily equivalent. This is a contradiction because in both cases the 
first algebra is totally atomic but the second one is not.}
\end{example}

{\em Acknowledgement:} I would like to express appreciation to
Prof. A. Katavolos for his helpful comments and suggestions during
the preparation of this work.
This research was partly supported by Special Account Research
Grant No. 70/3/7463 of the University of Athens.

A preliminary version of these results was presented by A. Katavolos
in the Operator Algebra Workshop held at
Queen's University Belfast in May 2006.

\end{document}